\documentclass[12pt]{article}

\usepackage{amsmath}
\usepackage{amssymb}
\usepackage[mathscr]{eucal}
\usepackage{graphicx}
\usepackage{color}

\newtheorem{Theorem}{\sc Theorem}
\newtheorem{Definition}[Theorem]{\sc Definition}
\newtheorem{Proposition}[Theorem]{\sc Proposition}
\newtheorem{Lemma}[Theorem]{\sc Lemma}
\newtheorem{Corollary}[Theorem]{\sc Corollary}

\newcommand{\cS}{\mbox{{${\cal S}$}}}
\newcommand{\R}{{\if mm {\rm I}\mkern -3mu{\rm R}\else \leavevmode
\hbox{I}\kern -.17em\hbox{R} \fi}}

\newcommand{\cF}{\mbox{{${\cal F}$}}}

\newcommand{\bu}{\mbox{\boldmath{$u$}}}
\newcommand{\bv}{\mbox{\boldmath{$v$}}}

\newcommand{\bx}{\mbox{\boldmath{$x$}}}

\newcommand{\fb}{\mbox{\boldmath{$f$}}}

\newcommand{\bvarepsilon}{\mbox{\boldmath{$\varepsilon$}}}
\newcommand{\bnu}{\mbox{\boldmath{$\nu$}}}

\newcommand{\bzero}{\mbox{\boldmath{$0$}}}

\def\sqr#1#2{{
    \vcenter{
         \vbox{\hrule height.#2pt
               \hbox{\vrule width.#2pt height#1pt \kern#1pt
                     \vrule width.#2pt
               }
               \hrule height.#2pt
         }
    }
}}

\def\real{\mathbb{R}}

\def\lista#1
{{ \itemindent 0.0cm \labelsep .2cm \leftmargin 0.8cm \rightmargin
0.0cm \labelwidth 0.6cm \topsep 0.0mm
\parsep 0.0mm
\itemsep 0.0mm
\begin{list}{}
{ \setlength{\leftmargin}{.8cm} \setlength{\rightmargin}{0.0cm}
\setlength{\parsep}{0.0mm} \setlength{\topsep}{.0mm}
\setlength{\parskip}{.0cm} \setlength{\itemsep}{.0cm} }
{#1}\end{list}} }

\newcounter{theorem}

\begin{document}

\title{\bf Tykhonov Well-posedness of  Elliptic Variational-Hemivariational Inequalities}

\vskip 16mm
\author{Mircea Sofonea$^{1,2}$\footnote{Corresponding author, E-mail : sofonea@univ-perp.fr}\, and\, Yi-bin Xiao$^{1}$
\\[6mm]		
{\it\small $^1$ School of Mathematical Sciences}\\
{\it \small University of Electronic Science and Technology of China}\\
{\it \small Chengdu, Sichuan, 611731, PR China}
\\[6mm]	
{\it \small $^2$ Laboratoire de Math\'ematiques et Physique}\\
{\it \small University of Perpignan Via Domitia}
\\
{\it\small 52 Avenue Paul Alduy, 66860 Perpignan, France}}

%



\date{}
\maketitle
\thispagestyle{empty}

\vskip 6mm

\noindent {\small{\bf Abstract.} \
We consider a class of elliptic variational-hemivaria\-tional inequalities in a abstract Banach space for which we  introduce the concept of well-posedness in the sense of Tykhonov.
We characterize the well-posedness in terms of metric properties of a family of associated sets.
Our results, which provide necessary and sufficient conditions for the well-posedness of inequalities under consideration, are valid under mild assumptions on the data.
Their proofs are based on arguments of monotonicity, lower semicontinuity and properties of the Clarke directional derivative. For well-posed inequalities we also prove a continuous dependence result of the solution with respect to the data. We illustrate our abstract results in  the study of one-dimensional examples, then we focus on some relevant particular cases, including  variational-hemivariational inequalities with strongly monotone operators. Finally, we consider a  model variational-hemivariational inequality  which arises in Contact Mechanics for which we discuss its well-posedness and provide the corresponding mechanical interpretations.}

\vskip 6mm

\noindent {\bf Keywords:} Tykhonov well-posedness, variational-hemivariational inequality, approximating sequence, contact problem, unilateral constraint, friction.

\vskip 6mm

\noindent {\bf 2010 Mathematics Subject Classification:} \ 49J40, 47J20, 49J45, 35M86, 74M10, 74M15.\\

\vskip 15mm

\section{Introduction}\label{s1}
\setcounter{equation}0


Well-posedness of  mathematical problems represents an important topic which was widely studied in the literature.
The concepts of well-posedness vary from problem to problem and  from author to author.  A few examples are the concept of well-posedness in the sense of Hadamard for partial differential equations, the concept of well-posedness in the sense of Tykhonov for minimization problems, the concept of well-posednes in the sense of Levitin-Polyak for constrained  optimization problems, among others.
Most of these concepts have been generalized in the recent years to various mathematical problems like inequality problems,  inclusion problems, fixed point problems, equilibrium problems and saddle point problems. The literature in the field is extensive, see for instance \cite {DZ,FHY,HYZ,L,PPXY,WXWC} and the references therein.

Inequality problems arise in the study of a large variety
of mathematical models used in Mechanics, Physics, Economy and Engineering Sciences. These models are usually expressed in terms of strongly nonlinear boundary value problems which, in a weak formulation, lead to
variational and hemivariational inequalities.
The theory of variational inequalities was developed in early sixty's, by using
arguments of monotonicity and convexity, including properties of the subdifferential 
of a convex function. References in the field 
include \cite{BC, ET, Gl, GLT,LXHC,LXH,SX2}.
The theory of hemivariational inequalities  started in early eighty's and used as main
ingredient the properties 
of the subdifferential in the sense of Clarke, defined for locally Lipschitz functions, 
which may be nonconvex. Comprehensive references on the field are the books \cite{MOSBOOK, NP,P}. 
Variational-hemivariational inequalities represent 
a special class of inequalities, in which both convex and nonconvex functions are involved. Recent references in the field is the monograph \cite{SMBOOK} as well as the papers \cite{Han,HS-ACTA,HSD,SXC1}. 

The present paper is devoted to the study of well-posedness in the sense of  Tykhonov for a class of elliptic variational-hemivariational inequalities.  Tykhonov well-posedness concept, introduced in \cite{Ty} for a minimization problem, is based on two main ingredients: the existence and uniqueness of the solution and the convergence to the unique solution of any approximating sequence. This concept was generalized to variational inequalities in \cite{LP1,LP2} and to hemivariational inequalities in \cite{GM}. References in the field  include \cite{SRX2,SRX1,XHW}. In particular, \cite{XHW} deals with the metric characterization of well-posedness of unconstrained hemivariational inequalities and inclusions, and
\cite{SRX2,SRX1} extend the results obtained there
to a special system of inequalities, the so-called split hemivariational and variational-hemivariational inequalities, respectively.

The present paper represents a continuation of  \cite{XHW} and  parallels  \cite{SRX2,SRX1}.  Thus, in contrast with  \cite{XHW}, here we deal with a more general type of inequalities which can be formulated as follows:  find an element  $u\in K$ such that
\begin{equation}\label{1}
\langle A u, v - u \rangle + \varphi (u, v) - \varphi (u, u) + j^0(u; v - u)
\ge \langle f, v - u \rangle \qquad\forall\,  v \in K.
\end{equation}

\medskip \noindent
In (\ref{1}) and everywhere in this paper, unless stated otherwise, $X$ is a real Banach space,   $\langle\cdot,\cdot\rangle$  denotes the duality pairing
between $X$ and its dual $X^*$,  $K$ is a nonempty subset of $X$,
$A \colon X \to X^*$,  $\varphi \colon X \times X \to \real$,
$j \colon X \to \real$ and $f\in X^*$. Note that the function
$\varphi(u, \cdot)$ is assumed to be convex
and the function $j$ is locally Lipschitz and, in general, is nonconvex. Moreover, $j^0(u;v)$ represents the general directional derivative of $j$ at the point $u$ in the direction $v$.  The inequality studied in \cite{XHW} represents a particular case of (\ref{1}), obtained when $K=X$ and $\varphi\equiv0$.



The existence and uniqueness of the solution of (\ref{1}) was proved  \cite{MOS30}, based
on arguments of  multivalued pseudomonotone operators
and the Banach fixed point theorem. The continuous dependence
of the solution with respect to the data  $A$, $\varphi$, $j$, $f$ and $K$ has been studied in \cite{XS3,ZLM}, where convergence results have been obtained, under various assumptions.
A comprehensive reference on the numerical analysis of  (\ref{1}) is the survey paper \cite{HS-ACTA}.
Nevertheless, the assumptions on the data considered in all these papers are quite strong. For instance, it is assumed that the space $X$ is reflexive, the operator $A$ is strongly monotone, $j$ satisfies the so-called relaxed monotonicity condition and, moreover, a smallness condition which relates $A$, $\varphi$ and $j$ is imposed. These assumptions represent sufficient conditions which guarantee the unique solvability of (\ref{1}).

Our aim in this paper is to study the well-posedness of the inequality (\ref{1}) in the sense of Tykhonov that we refer in what follows as well-posedness, for short.  We start by introducing this concept, then we provide necessary and sufficient
condition for the well-posedness. These conditions are expressed in terms of metric characterization of an useful set associated to the variational-hemivariational inequality. In comparison with  \cite{MOS30} (where only sufficient condition for the unique solvability of (\ref{1}) were considered), in this current paper we present a necessary and sufficient condition for its well-posedness (which implies its unique solvability), under less restrictive assumptions. This represents the first trait of novelty of our paper. In comparison with \cite{SRX2,SRX1,XHW} (where  particular inequalities of the form (\ref{1}) are considered),
in this paper we consider variational-hemivariational inequalities with constraints in which the function $\varphi$ depends on the solution  and prove that, under appropriate conditions, the well-posedness implies the continuous dependence of the solution with respect the data. This represents the second trait of novelty of our paper. Finally,  we illustrate  our results in the study of relevant particular cases and examples,
including an example which arises in Contact Mechanics.
This represents the third trait of novelty of our current paper.

The rest of the paper is structured as follows. In Section \ref{s2} we introduce the notion of well-posedness  for the variational-hemi\-variational inequality (\ref{1}), then we recall some preliminary material on nonsmooth analysis we need in the rest of the paper.
In Section \ref{s3} we state and prove our main results, Theorems \ref{t1} and \ref{t3}. Theorem  \ref{t1}  
provides necessary and sufficient conditions for the
well-posedness of inequality (\ref{1}). Theorem \ref{t3} provides a continuous dependence result of the solution with respect to the data. In Section \ref{s4} we provide two one-dimensional examples and, in Section \ref{s5} we consider two particular cases for which we specify our abstract results. The first one concerns a nonlinear equation and the second one a variational-hemivariational inequality.
Finally, in Section \ref{s6} we 
illustrate the use of our abstract results in the study of a model variational-hemivariational inequality 
which arises in Contact Mechanics and provide the corresponding mechanical interpretations.

\section{Problem statement and preliminaries}\label{s2}

In this section we introduce the problem statement and, to this end, we start with the following definitions.

\begin{Definition}\label{d7}
	A sequence $\{u_n\}\subset K$ is called an approximating sequence for the variational-hemivariational inequality $(\ref{1})$ if there exists a sequence  $\{\varepsilon_n\}\subset \real$ such that $0<\varepsilon_n\to 0$ and, for each $n\in\mathbb{N}$, the following inequality holds:
	\begin{eqnarray}
	&&\label{2}
	\langle A u_n, v - u_n \rangle + \varphi (u_n, v) - \varphi (u_n, u_n) + j^0(u_n; v - u_n)\\ [2mm]
	&&\qquad\ge \langle f, v - u_n \rangle-\varepsilon_n\| u_n - v \|_X
	\ \quad\forall\,  v \in K.\nonumber
	\end{eqnarray}
\end{Definition}

\begin{Definition}\label{d8} The variational-hemivariational inequality $(\ref{1})$ is said to be
{\rm well-posed} if it has a unique solution and every approximating sequence for  $(\ref{1})$ converges (strongly) in $X$ to the  solution.
\end{Definition}

Note that the concept of  well-posedness above extends that used 
in \cite{XHW} for pure hemivariational inequalities, but is quite different from that introduced in \cite{GM} for hemivariational inequalities with constraints.

\medskip

Our aim in what follows is to provide necessary and sufficient conditions which guarantee the well-posedness  of the variational-hemi\-variational inequality $(\ref{1})$. To this end, for each $\varepsilon>0$ we consider the set $\Omega(\varepsilon)$ defined as follows:
\begin{eqnarray}
&&\label{3}
\Omega(\varepsilon)=\{\, u\in K\ : \langle A u-f, v - u \rangle + \varphi (u, v) - \varphi (u, u) + j^0(u; v - u)\\ [2mm]
&&\hspace{38mm} \ge -\varepsilon\| u - v \|_X
\ \quad\forall\,  v \in K\,\}\nonumber
\end{eqnarray}
Moreover, we denote by ${\cal S}$ the set of solutions of inequality (\ref{1}), i.e.,
\begin{eqnarray}
&&\label{4}
\hspace{-15mm}{\cal S}=\{ u\in K : \langle A u-f, v - u \rangle + \varphi (u, v) - \varphi (u, u) + j^0(u; v - u)\ge0
\ \forall\,  v \in K\},
\end{eqnarray}
and we recall that ${\cal S}$ is said to be a singleton
if ${\cal S}$ has a unique element. Note that for each $\varepsilon>0$ the following inclusion holds:
\begin{equation}\label{5}
{\cal S}\subset\Omega(\varepsilon).
\end{equation}
Note that,  in general, this inclusion is strict.
The metric properties of the set $\Omega(\varepsilon)$ as $\varepsilon\to 0$  will play a crucial role in the Theorem \ref{t1}  we state and prove in the next section.  Its proof  requires some preliminaries of nonsmooth analysis  that we present in the rest of this section.
Everywhere in this paper
we use  $\|\cdot\|_X$ and $0_X$ for the norm and the zero element of space $X$, respectively. All the limits, upper and lower limits below are considered as $n\to\infty$, even if we do not mention it explicitly. The symbols ``$\rightharpoonup$"  and ``$\to$"
denote the weak and the strong convergence in the space $X$. 

\medskip
We start with some definitions related to the operator $A$ and functions $\varphi$, $j$.

\begin{Definition}\label{d1}
	An operator $A \colon X \to X^*$ is said to be:
	
	{\rm a) } monotone,
	if for all $u$, $v \in X$, we have $\langle Au - A v, u-v \rangle \ge 0$;

   {\rm	b)} strongly monotone, if there exists $m_A > 0$   such that
	\begin{equation}\label{sm}\langle Av_1 - Av_2, v_1 - v_2 \rangle
	\ge
	m_A \| v_1 - v_2 \|_X^{2} \qquad \forall\, v_1, v_2 \in X;
	\end{equation}
	
{\rm	c)} bounded, if $A$ maps bounded sets of $X$
	into bounded sets of $X^*$;
	
{\rm	d)} pseudomonotone,
	if it is bounded and $u_n \rightharpoonup u$  in $X$ with
	$$\displaystyle \limsup\,\langle A u_n, u_n -u \rangle \le 0$$
	implies\ \ 
	$$\displaystyle \liminf\, \langle A u_n, u_n - v \rangle\ge \langle A u, u - v \rangle\ \ \mbox{\rm for all}\ \, v \in X;$$
	
	
\end{Definition}

\begin{Definition}\label{d2}
	A function $\varphi \colon X \to \real$ is
	lower semicontinuous (l.s.c.) if $u_n \to u$ in $X$ implies 
	$\liminf \varphi (u_n)\ge\varphi(u)$.  A function $\varphi \colon X \to \real$ is
	weakly lower semicontinuous (weakly l.s.c.) if $u_n \rightharpoonup u$ in $X$ implies 
	$\liminf \varphi (u_n)\ge\varphi(u)$. 
\end{Definition}

\begin{Definition}\label{d3}
	 	A function $j \colon X \to \real$ is said to be 
		locally Lipschitz, 
		if for every 
		$u \in X$, there exists $N_u$ a neighborhood of $u$ and a constant $L_u>0$
		such that\ 
		$$
		|j(x) - j(y)| \le L_x \| x - y \|_X\qquad\forall\, x,\, y\in L_u.
		$$
	Assume in what follows that $j \colon X \to \real$ is a locally Lipschits function. Then, the
	generalized (Clarke) directional derivative of $j$ at the point
	$u \in X$ in the direction $v \in X$ is defined
	by
	\begin{equation*}
	j^{0}(u; v) = \limsup_{x \to u, \ \lambda \downarrow 0}
	\frac{j(x + \lambda v) - j(x)}{\lambda}.
	\end{equation*}
	The  generalized (Clarke) gradient (subdifferential) of $j$ at $u$
	is a subset of the dual space $X^*$ given by
	\begin{equation*}
	\partial j (u) = \{\, \xi \in X^* \mid j^{0}(u; v) \ge
	{\langle \xi, v \rangle} \quad \forall\, v \in X \, \}.
	\end{equation*}
\end{Definition}

\begin{Definition}\label{d4}	
	Let  $j:X\to\real$ be a locally Lipschitz functions.
	Then :
	
	{\rm a)\ } $j$ is said to be {\rm regular} (in the sense
	of Clarke) at the point $u \in X$ if for all $v \in X$ the one-sided directional
	derivative $j' (u; v)$ exists and $j^0(u; v) = j'(u; v)$.
	
	{\rm b)\ } $j$ is said to satisfy the relaxed monotonicity condition if there exists $\alpha_j>0$ such that
	\begin{equation*}\label{RR}
		\langle \xi_1 - \xi_2, u_1 - u_2 \rangle
		\ge -\alpha_j \, \| u_1 - u_2 \|_X^2 \qquad\forall\, u_i\in X,\ \xi_i\in\partial j(u_i),\ i=1,2.
	\end{equation*}
\end{Definition}

For the generalized (Clarke) directional derivative and the generalized (Clarke) gradient, we have the following properties, which could be found in \cite{Clarke, MOSBOOK}, for instance.

\begin{Proposition}\label{subdiff}
	Assume that $j \colon X \to \real$ is a locally Lipschitz function. Then the following hold:
	
	\medskip{\rm a)} 
	For every $u \in X$, the function
	$X \ni v \mapsto j^0(u;v) \in \real$ is positively  homogeneous and subadditive, 
	i.e., $j^0(u; \lambda v) = \lambda j^0(u; v)$ for all
	$\lambda \ge 0$, $v\in X$ and $j^0 (u; v_1 + v_2) \le
	j^0(u; v_1) + j^0(u; v_2)$ for all $v_1$, $v_2 \in
	X$, respectively.
	
	\medskip{\rm b)}
	For every $u,\ v \in X$, we have
	$j^0(u; v) = \max \, \{ \, \langle \xi, v \rangle
	\mid \xi \in \partial j(u) \, \}$.
	
	\medskip
	{\rm c)} The function $X \times X \ni (u, v)\mapsto j^0(u; v) \in \real$ is upper 
	semi-continuous, i.e., for all $u$, $v \in X$, $\{ u_n \}$, $\{ v_n \}
	\subset X$ such that $u_n \to u$ and $v_n \to v$ in $X$, we have $\limsup j^0(u_n; v_n) \le j^0(u; v)$.	
\end{Proposition}

\medskip
We also recall the following definition related to the metrics of the subsets in $X$.
\begin{Definition}\label{d5}
	Let $\Omega$ be a nonempty subset of $ X$. 
	Then, the diameter of $\Omega$, denoted ${\rm diam}(\Omega)$, is  defined by equality 
\[{\rm diam}(\Omega)= \displaystyle\sup_{a,\,b\in \Omega} \| a - b \|_X.\]
\end{Definition}

We we shall use this definition for  the set $\Omega(\varepsilon)$ defined by (\ref{3}).

\medskip





\section{Main results}\label{s3}

We  consider the following assumptions.
\begin{eqnarray}
&&\label{KK} K\ \mbox{is a closed subset of $X$.}\\ [3mm]
&&\label{AA} A:X\to X^*\ \mbox{is a pseudomonotone operator.}\\ [3mm]
&&\label{fifi} 
\left\{\begin{array}{ll} \varphi:X\times X\to\real\ \mbox{ and for all sequence}\ \{u_n\}\  \\[1mm]
\mbox{such that}\ \ u_n\to u\ \mbox{in}\ X\  \mbox{and all $v\in X$ we have} \\ [1mm]
\limsup\,\big(\varphi(u_n,v)-\varphi(u_n,u_n)\big)\le \varphi(u,v)-\varphi(u,u).
\end{array}\right.
\\ [3mm]
&&\label{jj} j:X\to \real\ \mbox{is a locally Lipschitz function.}
\end{eqnarray}

\medskip
Our first results in this section is  the following.

\medskip
\begin{Theorem}\label{t1}
	Let $X$ be a Banach space, $K$ a nonempty subset of $X$,
	$A \colon X \to X^*$,  $\varphi \colon X \times X \to \real$,
	$j \colon X \to \real$ and $f\in X^*$.  The following statements hold.

	\medskip
	{\rm a)} Under assumption  $(\ref{jj})$, the variational-hemivariational  inequality $(\ref{1})$ is well-posed if and only if its set of solution ${\cal S}$ is nonempty and\ {\rm diam}$(\Omega(\varepsilon))\to 0$\ as\ $\varepsilon\to 0$.

	\medskip
	{\rm b)}  Under assumptions  $(\ref{KK})$--$(\ref{jj})$, the variational-hemivariational  inequality $(\ref{1})$ is well-posed if and only if the set\ $\Omega(\varepsilon)$ is nonempty for each $\varepsilon>0$ and\ {\rm diam}$(\Omega(\varepsilon))$ $\to 0$\ as\ $\varepsilon\to 0$.

\end{Theorem}

\medskip\noindent{\it Proof.} a) Assume that  $(\ref{1})$ is well-posed. Then, by definition, ${\cal S}$ is  a singleton and, therefore, ${\cal S}\ne\emptyset$. Arguing by contradiction, we assume that \ {\rm diam}$(\Omega(\varepsilon))\not\to 0$\ as\ $\varepsilon\to 0$. Then, there exists $\delta_0\ge 0$, a sequence  $\{\varepsilon_n\}\subset\real$ and two sequences
$\{u_n\},\, \{v_n\}\subset X$ such that
 $0\ <\varepsilon_n\to0$, $u_n,\, v_n\in\Omega(\varepsilon_n)$ and
\begin{equation}\label{w1}
\|u_n-v_n\|_X\ge\frac{\delta_0}{2}\qquad\forall\, n\in\mathbb{N}.
\end{equation}
Now, since both $\{u_n\}$ and $\{v_n\}$ are approximating sequences for the variational-hemi\-variational inequality $(\ref{1})$, the well-posedness of $(\ref{1})$
implies that $u_n\to u$ and $v_n\to u$ in $X$ where $u$ denotes the unique element of ${\cal S}$. This is in contradiction with  $(\ref{w1})$. We conclude from here that {\rm diam}$(\Omega(\varepsilon))\to 0$\ as\ $\varepsilon\to 0$.

Conversely, assume that ${\cal S}$ is nonempty and\ {\rm diam}$(\Omega(\varepsilon))\to 0$\ as\ $\varepsilon\to 0$.
We claim that ${\cal S}$ is a singleton. Indeed, let $u,\, u'\in{\cal S}$ and let $\{u_n\}$ be an approximating sequence for $(\ref{1})$. Then there exists a sequence $\{\varepsilon_n\}\subset\real$ such that $0< \varepsilon_n\to0$ and $u_n\in\Omega(\varepsilon_n)$ for all $n\in\mathbb{N}$. We have
\[\|u-u'\|_X\le \|u-u_n\|_X+\|u'-u_n\|_X\le 2\,{\rm diam}(\Omega(\varepsilon_n))\to 0,\]
which implies that $u=u'$ and thus the claim is proved. Moreover, for any  approximating sequence we have
\[\|u-u_n\|_X\le \,{\rm diam}(\Omega(\varepsilon_n))\to 0,\]
which implies that $u_n\to u$ in $X$ and, therefore, (\ref{1}) is well-posed

\medskip
b) Assume that (\ref{1}) is well-posed. Then, we use
the part a) of the theorem and inclusion (\ref{5}) to see that the set\ $\Omega(\varepsilon)$ is nonempty for each $\varepsilon>0$ and\ {\rm diam}$(\Omega(\varepsilon))\to 0$\ as\ $\varepsilon\to 0$.

Conversely, assume that the set\ $\Omega(\varepsilon)$ is nonempty for each $\varepsilon>0$ and\ {\rm diam}$(\Omega(\varepsilon))\to 0$\ as\ $\varepsilon\to 0$.
 Then, using (\ref{5}), again,  and Definition \ref{d5} we deduce that
 \begin{equation}\label{xx}
  {\cal S}\ne\emptyset\quad\Longrightarrow\quad {\cal S}\quad\mbox{is a singleton.}
 \end{equation}
 We now prove the solvability of the variational-hemivariational inequality (\ref{1}) and, to this end, we use a pseudomonotonicity argument. Let
 $\{u_n\}$ be an approximating sequence for $(\ref{1})$. Then there exists a sequence $\{\varepsilon_n\}\subset\real$ such that $0< \varepsilon_n\to0$ and $u_n\in\Omega(\varepsilon_n)$ for all $n\in\mathbb{N}$. Since \ {\rm diam}$(\Omega(\varepsilon_n))\to 0$
 it follows that
 $\{u_n\}$ is a Cauchy sequence in $X$ and, therefore, there exists $u\in X$ such that
 \begin{equation}\label{w2}
 u_n\to u\quad{\rm in}\quad X.
 \end{equation}
 This convergence combined with assumption (\ref{KK}) yields
 \begin{equation}\label{w3}
 u\in K.
 \end{equation}
 We now use (\ref{2}) to see that
 \begin{eqnarray*}
 &&
 \langle A u_n, u_n-v \rangle \le \varphi (u_n, v) - \varphi (u_n, u_n) + j^0(u_n; v - u_n)\\ [2mm]
 &&\qquad+\langle f, u_n-v \rangle+\varepsilon_n\| u_n - v \|_X
 \ \quad\forall\,  v \in K,\ n\in\mathbb{N}.\nonumber
 \end{eqnarray*}
 Next, we pass to the upper limit as $n\to\infty$ in this inequality and use the convergence (\ref{w2}), assumption (\ref{fifi}) and Proposition \ref{subdiff} c) to deduce that
 \begin{eqnarray}
 \label{w4}\limsup\, \langle A u_n, u_n-v \rangle
  &\le& \varphi (u, v) - \varphi (u, u)+j^0(u; v- u)\\ [2mm] 
 &&-\langle f, v-u \rangle
 \qquad\forall\,  v \in K.\nonumber
 \end{eqnarray}
 
 On the other hand, regularity (\ref{w3}) allows us to test with $v=u$ in (\ref{w4}) to find that
 \[\limsup\, \langle A u_n, u_n-u \rangle\le 0.\]
Therefore, by the pseudomonotonicity of the operator $A$, guaranteed by  assumption (\ref{AA}), we obtain
\begin{equation}\label{w5}
\liminf\, \langle A u_n, u_n - v \rangle\ge \langle A u, u - v \rangle\qquad\forall\, v\in X.
\end{equation}
We now combine (\ref{w3}), (\ref{w4}) and (\ref{w5}) to see that $u$ is a solution to the variational-hemivariational inequality (\ref{1}), i.e., $u\in\cS$.  We now use (\ref{xx}) and  (\ref{w2}) to see that ${\cal S}$ is a singleton and any approximating sequence of (\ref{1}) converges to the unique element of ${\cal S}$.  It follows from here that the variational-hemivariational  inequality $(\ref{1})$ is well-posed, which concludes the proof.
\hfill$\Box$

\medskip

Consider  now the sequences  $\{\varphi_n\}$,  $\{j_n\}$, $\{f_n\}$ such that	that, for each $n\in\mathbb{N}$,  the following conditions hold:
\begin{eqnarray}
&&\hspace{-9mm}\label{fin}\left\{\begin{array}{ll} \varphi_n\colon X\times X \to\real\ \mbox{and there exists $b_n>0$  such that}\\[2mm]
\varphi_n(u,v)-\varphi_n(u,u)-\varphi(u,v)+\varphi(u,u)\le  b_n\|u-v\|_X \ \ \forall\, u,\, v\in X.
\end{array}\right.\\[3mm]
&&\hspace{-9mm}\label{jn}\left\{\begin{array}{ll} j_n\colon X \to\real\ \mbox{ is a locally Lipschitz function and}\\ [2mm] 
\mbox{there exists $c_n>0$  such that}\\[2mm]
j_n^0(u;v-u)-j^0(u,v-u)\le c_n\|u-v\|_X \ \ \forall\, u,\, v\in X.
\end{array}\right.
\\[3mm]
&&\label{fn}f_n\in X^*.
\end{eqnarray}

With these data, for each $n\in\mathbb{N}$, we consider the following problem : find $u_n\in K$ such that
\begin{eqnarray}
&&\label{1n}
\langle Au_n, v - u_n \rangle + \varphi_n(u_n, v) - \varphi_n(u_n, u_n) + j^0_n(u_n; v - u_n)\\ [2mm]
&&\qquad\qquad\ge\langle f_n, v - u_n \rangle
\ \quad\forall\,  v \in K.\nonumber
\end{eqnarray}
Finally, we  assume that
\begin{eqnarray}
&&\label{cn}b_n\to 0,\\[2mm]
&&\label{ccn}c_n\to 0,\\[2mm]
&&\label{cfn}f_n\to f \quad{\rm in}\quad X^*.
\end{eqnarray}

Then, we have the following result.

\begin{Theorem}\label{t3}
Let $X$ be a Banach space,  $K$ a nonempty subset of $X$,
$A \colon X \to X^*$,  $\varphi \colon X \times X \to \real$,
$j \colon X \to \real$ and $f\in X^*$. Assume that $(\ref{jj})$, $(\ref{fin})$--$(\ref{fn})$, $(\ref{cn})$--$(\ref{cfn})$ hold and the variational-hemivariational inequality $(\ref{1})$ is  well-posed. Moreover, for each $n\in\mathbb{N}$, let $u_n$ be a solution of inequality  $(\ref{1n})$.
Then $u_n\to u$ in $X$. 
\end{Theorem}

\noindent{\it Proof.} Let $n\in\mathbb{N}$ and $v\in X$. We write
\begin{eqnarray*}
&&
\langle Au_n, v - u_n \rangle + \varphi(u_n, v) - \varphi(u_n, u_n) + j^0(u_n; v - u_n)-\langle f, v - u_n \rangle\\  [2mm]
&&= \langle Au_n, v - u_n \rangle+
\\ [2mm]
&&+\varphi_n(u_n, v) - \varphi_n(u_n, u_n)+ \varphi_n(u_n, u_n) - \varphi_n(u_n,v)+\varphi(u_n, v) - \varphi(u_n, u_n)
\\ [2mm]
&&+j^0_n(u_n; v - u_n)-j^0_n(u_n; v - u_n)+j^0(u_n; v - u_n)
\\ [2mm]
&&-\langle f_n, v - u_n \rangle+ \langle f_n, v - u_n \rangle -\langle f, v - u_n \rangle,
\end{eqnarray*}
then we use assumptions  (\ref{fin}) and (\ref{jn}) to deduce that
\begin{eqnarray*}
&&\langle Au_n, v - u_n \rangle + \varphi(u_n, v) - \varphi(u_n, u_n) + j^0(u_n; v - u_n)-\langle f, v - u_n \rangle\\ [2mm]
&&\ge\langle Au_n, v - u_n \rangle+\varphi_n(u_n, v)- \varphi_n(u_n, u_n)+j^0_n(u_n; v - u_n)-\langle f_n, v - u_n \rangle
\\ [2mm]
&&-b_n\|u_n-v\|_X -c_n\|u_n-v\|_X -\|f_n-f\|_{X^*} \|u_n-v\|_X. 
\end{eqnarray*}
Moreover, exploiting (\ref{1n}) we find that
\begin{eqnarray}
	&&\label{z1}\langle Au_n, v - u_n \rangle + \varphi(u_n, v) - \varphi(u_n, u_n) + j^0(u_n; v - u_n)-\langle f, v - u_n \rangle\\ [2mm]
	&&\qquad\qquad\qquad\ge-\varepsilon_n \|u_n-v\|_X \nonumber
\end{eqnarray}
where
\[\varepsilon_n=b_n+c_n+\|f_n-f\|_{X^*}.\]
Note that assumptions $(\ref{cn})$--$(\ref{cfn})$ imply that
$\varepsilon_n\to 0$ and, therefore, inequality (\ref{z1}) and  Definition \ref{d7} show that $\{u_n\}$ is an approximating sequence for the variational-hemivariational inequality  (\ref{1}).  Theorem \ref{t3} is now a direct consequence of Definition \ref{d8}.
\hfill$\Box$

\medskip
In the next result, for any $f\in X^*$, we denote by
$u(f)$ the solution of the variational-hemivatriational inequality (\ref{1}), assumed to be unique.

\begin{Corollary}\label{c1}
Let $X$ be a Banach space,  $K$ a nonempty subset of $X$
$A \colon X \to X^*$,  $\varphi \colon X \times X \to \real$,
$j \colon X \to \real$ and assume that $(\ref{jj})$ holds. Moreover, assume that the variational-hemivariational inequality $(\ref{1})$ is well-posed for any $f\in X^*$. Then the operator $f\mapsto u(f):X^*\to X$ is continuous.

\end{Corollary}

\medskip\noindent{\it Proof.}
Corollary \ref{c1} follows directly form Theorem \ref{t3} and  Definition \ref{d7}, by considering  the particular case when $\varphi_n=\varphi$
and $j_n=j$.\hfill$\Box$

\section{Two one-dimensional examples}\label{s4}

It follows from Theorem \ref{t1} that 
the well-posedness of the variational-hemivariational inequality (\ref{1})  is related to the properties of the sets $\Omega(\bvarepsilon)$ defined by (\ref{3}). Note  that, in general, it is not easy to describe explicitly these sets.
The two examples we present in what follows have the merit that in each case  we can clearly determinate the sets $\Omega(\bvarepsilon)$. This allows us to use  Theorem \ref{t1} in order to see that, for some data $f$, the corresponding inequalities are well-posed and, for other data, they fail to be.

\medskip\noindent {\bf a) First one-dimensional example.} The example is the following:
$X=K=\real$, $Au=u$, $\varphi(u)=0$ for all $u\in \real$ and
\begin{equation*}
j(u)=\left\{\begin{array}{ll} \frac{1}{2}u^2\hspace{25mm}{\rm if}\quad u<1,\\[2mm]
2u-\frac{1}{2}u^2-1\quad\quad{\rm if}\quad 1\le u\le 2,\\ [2mm]
\frac{1}{2}u^2-1\hspace{18mm}{\rm if}\quad u>2.
\end{array}\right.
\end{equation*}
It is easy to see that this function is locally Lipschitz yet nonconvex. Moreover, it is regular and a simple calculation shows that
$j^0(u;v)=p(u)v$ for all $u,\, v\in\real$ where $p:\real\to\real$ is the function defined by
\begin{equation}\label{pp}
p(u)=\left\{\begin{array}{ll} u\hspace{16mm}{\rm if}\quad u<1,\\[2mm]
2-u\quad\quad{\rm if}\quad 1\le u\le 2,\\ [2mm]
u\hspace{16mm}{\rm if}\quad u>2.\\[2mm]
\end{array}\right.
\end{equation}
Therefore, in this particular case
the variational-hemivariational inequality (\ref{1}) reads
\begin{equation}\label{2pp}
u(v-u)+p(u)(v-u)\ge f(v-u)\qquad\forall\,v\in\real
\end{equation}
or, equivalently,
\begin{equation}\label{nep}
u+p(u)=f.
\end{equation}
Now, using (\ref{pp}) we find  that the solution of
the nonlinear equation (\ref{nep}) is given by
\begin{equation}\label{sp}
u=\left\{\begin{array}{ll} \frac{f}{2}\hspace{14mm}{\rm if}\quad f<2,\\[2mm]
[1,2]\quad\quad{\rm if}\quad f=2,\\ [2mm]
\frac{f}{2}\hspace{14mm}{\rm if}\quad f>2.\\[2mm]
\end{array}\right.
\end{equation}
This shows that equation (\ref{2pp}) has a unique solution if $f<2$ or $f>2$ and for $f=2$ it has an infinity  of solutions, since in this case any element $u\in[1,2]$ is a solution to (\ref{2pp}).

We now specify  the set $\Omega(\varepsilon)$ for inequality (\ref{2pp}), for any $\varepsilon >0$. Using the arguments above and denoting $w=v-u$ we see that $u\in \Omega(\varepsilon)$ if and only if
\begin{equation}\label{ne}
(u+p(u)-f)w+\varepsilon|w|\ge 0\qquad\forall\,w\in \real.
\end{equation}
Next, based on the elementary equivalence
\[xw+\varepsilon|w|\ge 0\quad\forall\,w\in\real\quad\Longleftrightarrow\quad x\in[-\varepsilon,\varepsilon],\]
we deduce that (\ref{ne}) is equivalent to the inequality
\begin{equation}\label{npp}
-\varepsilon\le u+p(u)-f\le\varepsilon.
\end{equation}
We now use the (\ref{pp}) to see that
\begin{equation*}
u+p(u)-f=\left\{\begin{array}{ll} 2u-f\hspace{16mm}{\rm if}\quad u<1,\\[2mm]
2-f\hspace{18mm}{\rm if}\quad 1\le u\le 2,\\ [2mm]
2u-2-f\qquad{\rm if}\quad u>2.\\[2mm]
\end{array}\right.
\end{equation*}
Then, using a graphic method it is easy to see that
\begin{equation}\label{app}
\Omega(\varepsilon)=\left\{\begin{array}{ll} \Big[\frac{f-\bvarepsilon}{2},\frac{f+\bvarepsilon}{2}\Big]\qquad{\rm if}\quad f<2\quad{\rm and}\quad \varepsilon<2-f,\\[4mm]
\Big[\frac{2-\bvarepsilon}{2},\frac{\bvarepsilon+4}{2}
\Big]\qquad{\rm if}\quad f=2 \quad \forall\,\bvarepsilon>0,\\ [4mm]
\Big[\frac{f-\bvarepsilon}{2},\frac{f+\bvarepsilon}{2}
\Big]\qquad{\rm if}\quad f>2 \quad{\rm and}\quad \varepsilon<f-2.\\[2mm]
\end{array}\right.
\end{equation}
It follows from here that
\begin{equation}\label{bpp}
\lim_{\bvarepsilon\to 0}{\rm diam}(\Omega(\varepsilon))=\left\{\begin{array}{ll} 0\qquad{\rm if}\quad f<2,\\[4mm]
1\qquad{\rm if}\quad f=2 ,\\ [4mm]
0\qquad{\rm if}\quad f>2.
\end{array}\right.
\end{equation}
We now apply Theorem \ref{t1} to see that  the hemivariational inequality (\ref{2pp}) is well-posed
if and only if $f\ne 2$. This result is in agreement with 
our previous computations since, recall, the solution of inequality (\ref{2pp}) is given by (\ref{sp}).

\medskip\noindent {\bf b) Second one-dimensional example.} Our second example is the following:
$X=K=\real$, $Au=u$, $\varphi(u)=0$ for all $u\in \real$ and
\begin{equation*}
j(u)=\left\{\begin{array}{ll} -u^2+3u\qquad{\rm if}\quad u<1,\\[2mm]
u+1\hspace{16mm}{\rm if}\quad u\ge 1.
\end{array}\right.
\end{equation*}
This function is locally Lipschitz yet nonconvex. Moreover, it is regular and a simple calculation shows that
$j^0(u;v)=p(u)v$ for all $u,\, v\in\real$ where $p:\real\to\real$ is the function defined by
\begin{equation}\label{mpp}
p(u)=\left\{\begin{array}{ll} -2u+3\qquad{\rm if}\quad u<1,\\[2mm]
1\hspace{21mm}{\rm if}\quad u\ge 1.
\end{array}\right.
\end{equation}
Therefore, using the arguments in the previous example, it follows that in this case the set of solutions to the variational-hemivariational inequality 
(\ref{1}) is given by
\begin{equation}\label{msp}
{\cal S}=\left\{\begin{array}{ll}
\emptyset\hspace{31mm}{\rm if}\quad f<2,\\[2mm]
{1}\hspace{31mm}{\rm if}\quad f=2,\\ [2mm]
\{3-f,f-1\}\quad\quad{\rm if}\quad f>2.\\[2mm]
\end{array}\right.
\end{equation}
This shows that the inequality has a unique solution if $f=2$, two solutions if $f>2$, and no solution if $f<2$.

The set $\Omega(\varepsilon)$ can be determined, for any $\varepsilon>0$, by using arguments similar to those used in the previous example. We have
\begin{equation*}
\Omega(\varepsilon)=\left\{\begin{array}{ll} \emptyset\hspace{60mm}\qquad{\rm if}\ f<2\ {\rm and}\ \varepsilon<2-f,\\[4mm]
[1-\varepsilon,1+\varepsilon]\hspace{48mm}{\rm if}\ f=2 \ \ \forall\,\bvarepsilon>0,\\ [4mm]
[3-f-\bvarepsilon,3-f+\bvarepsilon
]\cup[f-1-\bvarepsilon,f-1+\bvarepsilon
]\\
\hspace{69mm}\ {\rm if}\ f>2 \ {\rm and}\ \varepsilon<f-2.
\end{array}\right.
\end{equation*}
It follows from here that \[\lim_{\bvarepsilon\to 0}{\rm diam}(\Omega(\varepsilon))=0\quad{\rm  if}\quad f=2\] 
and
\[\lim_{\bvarepsilon\to 0}{\rm diam}(\Omega(\varepsilon))=2f-4>0\quad{\rm  if}\quad f>2.\]

We now apply Theorem \ref{t1} to see that  the corresponding  inequality (\ref{1}) is well-posed
if and only if $f=2$. This result is in agreement with 
our previous computations since, recall, the solution of this inequality  is given by (\ref{msp}).  Moreover, it is in contrast with the situation in our previous example since inequality (\ref{2pp}) is well-posed  if and only if $f\ne2$.

\section{Two relevant particular cases}\label{s5}

In this section we present two relevant particular cases  of variational-hemivariational inequalities of the form (\ref{1}) for which we apply and complete the  results in Theorem \ref{t1}. The problems we consider here have some interest on their own. 

\medskip\noindent
{\bf a) A nonlinear equation in reflexive Banach spaces.} For the particular case we consider in this subsection that $X$ is a reflexive Banach space and the norm on $X$ is strictly convex. Note that this is not a restriction since it is well know that
any reflexive Banach space $X$
can be always considered as equivalently renormed strictly convex space. Moreover, we assume that $K=X$ and there exist two operators $L:X\to X^*$ and $P:X\to X^*$ such that $\varphi(u,v)=\langle Lu,v\rangle$, $j^0(u,v)=\langle Pu,v\rangle$ for all $u,\, v\in X$. Examples of such functions were already given in Section \ref{s4} and another example of $j$ with this property will be provided in Section \ref{s6}, too. Denote by $T:X\to X^*$ the operator given by $T=A+L+P$ .
Then it is easy to see that the variational-hemivariational inequality  (\ref{1}) is equivalent to the variational inequality
\begin{equation}\label{e1}
u\in X,\qquad\langle Tu,v-u\rangle\ge \langle f,v-u\rangle\qquad\forall\, v\in X,
\end{equation}
which, in turn, is equivalent with the equation
\begin{equation}\label{e2}
Tu=f.
\end{equation}
Based on this equivalence we transpose all the definitions and notions related to the well-posedness of inequality (\ref{e1}) to  corresponding  definitions and notions for
equation (\ref{e2}). For instance, 
we say that equation (\ref{e2}) is well-posed if the variational inequality (\ref{e1}) is well-posed in the sense of Definition \ref{d7}.  Moreover, using (\ref{3}) we see that for any $\varepsilon>0$ the set $\Omega(\varepsilon)$ associated to (\ref{e2}) is given by
\begin{eqnarray}
&&\label{3e}
\Omega(\varepsilon)=\{\, u\in X\ :\ \langle T u-f, w\rangle + \varepsilon\|w\|_X\ge 0
\ \quad\forall\,  w\in X\,\}.
\end{eqnarray}

For any $\theta\in X^*$  and $\varepsilon>0$ we denote in what follows by $\overline{B}(\theta,\varepsilon)$ the closed ball of center $\theta$ and radius $\varepsilon$ in the dual space $X^*$, i.e., 
\begin{eqnarray}
&&\label{3eb}
\overline{B}(\theta,\varepsilon)=\{\, \xi\in X^*\ :\ \|\xi-\theta\|_{X^*}\le\varepsilon\,\}.
\end{eqnarray}
Then, we have the following characterization of the the set $\Omega(\varepsilon)$ given by (\ref{3e}).

\begin{Proposition}\label{p8}
	Let $X$ be a reflexive strictly convex Banach space. Then, for each $\varepsilon>0$ the following equivalence holds:
	\begin{equation}\label{e5}
	u\in \Omega(\varepsilon) \quad\Longleftrightarrow\quad
	Tu\in \overline{B}(f,\varepsilon).
	\end{equation}
\end{Proposition}

\noindent{\it Proof.} We claim that for each $\varepsilon>0$ the following equivalence holds:
\begin{equation}\label{e6}
x^*\in X^*, \quad \langle x^*,w\rangle+\varepsilon\|w\|_X\ge 0\quad\forall\,w\in X \quad\Longleftrightarrow\quad
x^*\in \overline{B}(0_{X^*},\varepsilon).
\end{equation}

Indeed, assume that $x^*$ is an element of $X^*$ such that
\[\langle x^*,w\rangle+\varepsilon\|w\|_X\ge 0\quad \forall\,w\in X.\] By reflexivity of $X$ we know that there exists an element $\theta\in X$ such that
$\langle x^*,\theta\rangle=\|\theta\|_X^2=\|x^*\|_{X^*}^2$
and, testing with $ w=-\theta$ in the previous inequality we deduce that $\|x^*\|_{X^*}\le\varepsilon$, i.e., $x^*\in \overline{B}(0_{X^*},\varepsilon)$.
Conversely, if $x^*\in \overline{B}(0_{X^*},\varepsilon)$ we have  $\|x^*\|_{X^*}\le\varepsilon$ and, therefore, for each $w\in X$, we deduce that 
\[\langle x^*,w\rangle+\varepsilon\|w\|_X\ge-\|x^*\|_{X^*}\|w\|_X+\varepsilon\|w\|_X=(\varepsilon-\|x^*\|_{X^*})\|w\|_X\ge 0,\]
which concludes the proof of the claim.

Let $\varepsilon>0$.  We use definition (\ref{3e})
and equivalence (\ref{e6}) so see that
\begin{eqnarray*}
&&u\in \Omega(\varepsilon) \quad\Longleftrightarrow\quad
\langle T u-f, w\rangle + \varepsilon\|w\|_X\ge 0
\ \quad\forall\,  w\in X\\[2mm]
&& \qquad\Longleftrightarrow\quad Tu-f\in\overline{B}(0_{X^*},\varepsilon)\quad\Longleftrightarrow\quad Tu\in\overline{B}(f,\varepsilon),
\end{eqnarray*}
which concludes the proof.

\hfill$\Box$

We now use this result in order to give a characterization of the well-posedness of the equation  (\ref{e2}).

\begin{Theorem}\label{t9}
	Let $X$ be a reflexive Banach space and let $T:X\to X^*$. Then 
	equation  $(\ref{e2})$ is well-posed for each $f\in X^*$ if and only if the operator $T$ is invertible  and its inverse $T^{-1}:X^*\to X$ is continuous.
\end{Theorem}

\medskip\noindent{\it Proof.}  Assume that (\ref{e2})  is well-posed, for any $f\in X^*$. Then, it follows that for each $f\in X^*$ there exists a unique element $u\in X$ such $Tu=f$ and, therefore, the operator $T$ is invertible. Let $\{f_n\}\subset X^*$, $f\in X^*$ be such that $f_n\to f$ in $X^*$, let $u=T^{-1}f$ and, for each $n\in\mathbb{N}$, denote $u_n=T^{-1}f_n$, 
$\varepsilon_n=\|f_n-f\|_{X^*}$. We have
\[\|Tu_n-f\|_{X^*}=\|f_n-f\|_{X^*}=\varepsilon_n\]
and, therefore $Tu_n\in \overline{B}(f,\varepsilon)$ for each $n\in\mathbb{N}$. We now use  (\ref{e5}) to see that $u_n\in\Omega(\varepsilon_n)$, for each $n\in\mathbb{N}$.
On the other hand, the convergence $f_n\to f$ in $X^*$ guarantees that $\varepsilon_n\to 0$ which shows that $\{u_n\}$ is an approximating sequence for the equation (\ref{e2}). Since  by assumption this equation is well-posed and its solution is  $u$, we have $u_n\to u$ in $X$,  i.e., $T^{-1}f_n\to T^{-1}f$ in $X$.
We conclude from above that $T^{-1}:X^*\to X$ is a  continuous operator.

Conversely, assume that $T$ is invertible  and its inverse $T^{-1}:X^*\to X$ is continuous. Let $f\in X^*$ and let 
$T^{-1}f=u$ or, equivalently, $Tu=f$. Let  $\{u_n\}$ be an approximating sequence for (\ref{e2}). The, by definition, there exists a sequence $\{\varepsilon_n\}\subset\real$ such that  $0<\varepsilon_n\to 0$ and $u_n\in\Omega(\varepsilon_n)$, for each $n\in\mathbb{N}$. 
Note that equivalence (\ref{e5}) yields $Tu_n\in\overline{B}(f,\varepsilon_n)$ or, equivalently,
$\|Tu_n-f\|_{X^*}\le\varepsilon_n$ for each $n\in\mathbb{N}$, which shows that $Tu_n\to f$ in $X^*$. Using now the continuity of $T^{-1}$ we find that $T^{-1}(Tu_n)\to T^{-1} f$ in $X$ and, therefore, $u_n\to u$ in $X$. This shows that equation  $(\ref{e2})$ is  well-posed.
\hfill$\Box$

\medskip\noindent{\bf b) Variational-hemivariational inequalities with strongly monotone operators.}
We now  study the well-posedness of the variational-hemivariational inequality (\ref{1}) in the particular case when  the operator $A$ is strongly monotone. 
The complete list of assumptions we consider on the data  is the following.
\begin{eqnarray}
&&\label{K}
K \ \mbox{is nonempty, closed and convex subset of} \ X .
\\[2mm]
&&\label{A1}\left\{\begin{array}{l} A:X\to X^*\ \mbox{is a strongly monotone operator,}\\
\mbox{i.e., it satisfies condition (\ref{sm}) with}\ m_A>0.
\end{array}
\right.\\[2mm]
&&\label{A2} A:X\to X^*\ \mbox{is a pseudomonotone operator.}
\end{eqnarray}
\begin{eqnarray}
&&\left\{\begin{array}{l}
\varphi \  \colon X \times X \to \real \ \mbox{and there
exists} \ \alpha_\varphi > 0 \ \mbox{such that}  \\ [2mm]
\varphi (\eta_1, v_2) - \varphi (\eta_1, v_1) + \varphi (\eta_2, v_1) - \varphi (\eta_2, v_2)\\ [2mm]
\quad\le \alpha_\varphi \| \eta_1 - \eta_2 \|_X \, \| v_1 - v_2 \|_X\ \ 
\ \forall\, \eta_1, \eta_2, v_1, v_2 \in X.
\end{array}
\right.
\label{fi1}
\end{eqnarray}
\begin{eqnarray}
&&\varphi (\eta, \cdot) \colon X \to \real \ \mbox{is convex and l.s.c. for all}\ \eta \in X.  
\label{fi2}
\end{eqnarray}
\begin{equation}
\label{j1} \left\{
\begin{array}{l}
\mbox{there exists} \ \alpha_j > 0 \ \mbox{such that} \\ [2mm]
j^0(v_1; v_2 - v_1) + j^0(v_2; v_1 - v_2) \le \alpha_j \, \| v_1 - v_2 \|_X^2 \\ [2mm]
\quad \forall\,  v_1, v_2 \in X.
\end{array}
\right.
\end{equation}
\smallskip
\begin{equation}
\label{j2} \left\{
\begin{array}{l}
\mbox{there exists} \ \ c_0, c_1 \ge 0 \ \mbox{such that}\\[2mm] 
\| \xi \|_{X^*} \le c_0 + c_1 \, \| v \|_X \quad \forall\,
 v \in X,\ \xi\in \partial j(v).\hspace{7mm}
\end{array}
\right.
\end{equation}
\begin{eqnarray}
&&\label{smal}
\alpha_\varphi + \alpha_j < m_A.
\end{eqnarray}

\medskip

It can be proved that for a locally Lipschitz function $j \colon X \to \real$,
hypothesis~$(\ref{j1})$ is equivalent to the so-called relaxed monotonicity condition introduced in Definition \ref{d4}(b). A proof of the statement can be found in, e.g., {\rm \cite{MOSBOOK}}.
Note also that if $j \colon X \to \real$ is a convex function,
then $(\ref{j1})$ 
holds  with $\alpha_j = 0$, since it reduces
to the monotonicity of the (convex) subdifferential.
Examples of functions which satisfy conditions
(\ref{jj}), (\ref{j1}),  (\ref{j2}) have been provided in \cite{MOSBOOK,MOS30}, for instance.

\medskip
\begin{Theorem}\label{t4} Let $X$ be a Banach space,  $K$ a nonempty subset of $X$,
	$A \colon X \to X^*$,  $\varphi \colon X \times X \to \real$,
	$j \colon X \to \real$ and assume that $(\ref{jj})$, $(\ref{A1})$, $(\ref{fi1})$, $(\ref{j1})$ and $(\ref{smal})$ hold. Then, for all $f\in X^*$,  the following statements are equivalent:
	
\medskip	
{\rm  a)} The variational-hemivariational inequality $(\ref{1})$ has a unique solution.

\medskip
{\rm  b)}  The variational-hemivariational inequality $(\ref{1})$ is well-posed.
\end{Theorem}

\noindent
{\it Proof.} Assume a). Let $u\in K$ be the unique solution of (\ref{1}) and let $\{u_n\}\subset K$ be an approximating sequence. Let $n\in\mathbb{N}$. We write (\ref{1}) with $v=u_n$, (\ref{2}) with $v=u$, then we add the resulting inequalities to see that
\begin{eqnarray*}
&&\langle A u_n-Au, u_n -u\rangle \le \varphi (u_n, u) - \varphi (u_n, u_n) + \varphi (u, u_n)-\varphi (u, u)\\ [2mm]
&&\quad+j^0(u_n; u - u_n)+j^0(u; u_n - u)+\varepsilon_n\| u_n - u \|_X.
\end{eqnarray*}
We now use assumptions $(\ref{A1})$, $(\ref{fi1})$ and $(\ref{j1})$ to obtain that
\begin{equation*}
m_A\| u_n - u \|_X^2\le \alpha_\varphi\| u_n - u \|_X^2+ \alpha_j\| u_n - u \|_X^2+\varepsilon_n\| u_n - u \|_X.
\end{equation*}
Therefore, the smallness condition (\ref{smal}) yields
\begin{equation*}
\| u_n - u \|_X\le \frac{\varepsilon_n}{m_A-\alpha_\varphi-\alpha_j}
\end{equation*}
and, since Definition \ref{d7} guarantees that $\varepsilon_n\to 0$, we deduce that $u_n\to u$ in $X$.
This proves that  {\rm b)} holds. We conclude from here that {\rm a)} $\Longrightarrow$ {\rm b)}.  Note that
the converse implication, {\rm b)} $\Longrightarrow$ {\rm a)}, is a direct consequence of Definition \ref{d8}. We conclude from above that the statements a) and b) are equivalent, which completes the proof.
\hfill$\Box$

\medskip
Note that Theorem \ref{t4} provides an equivalence result. It does not guarantees that the statements {\rm a)}, {\rm b)} above are valid. Sufficient conditions  which guarantee the validity of these statements are provided by the following result.

\begin{Theorem}\label{t5}
	Assume  that $X$ is a reflexive Banach space and, moreover, assume that $(\ref{jj})$, $(\ref{K})$--$(\ref{smal})$ hold. Then, for each $f\in X^*$, the variational-hemivariational inequality~$(\ref{1})$ has a unique solution $u=u(f) \in K$.
\end{Theorem}

A proof of   Theorem~\ref{t5} can be found \cite{MOS30}, see also Remark 13 on \cite{SMBOOK}. It is carried out in several steps, by using the properties of the subdifferential, a surjectivity result for  pseudomonotone multivalued operators and the Banach fixed point argument.

\medskip
Using now  Theorems \ref{t4}, \ref{t5} and Corollary \ref{c1} it is easy to deduce the following result

\begin{Corollary}\label{c2}
Assume  that $X$ is a reflexive Banach space and, moreover, assume that $(\ref{jj})$, $(\ref{K})$--$(\ref{smal})$ hold. Then, for all $f\in X^*$, the variational-hemivariational inequality~$(\ref{1})$ is well-posed. Moreover, the operator $f\mapsto u(f):X^*\to X$ is continuous.
\end{Corollary}

We end this section with the the remark that, under assumptions in Corollary \ref{c2}, it can be proved that
the operator $f\mapsto u(f):X^*\to X$ is Lipschitz continuous.

\section{An application to Contact Mechanics}\label{s6}

The results presented in Sections \ref{s3} and \ref{s5} can be used in the study of various mathematical models which describe the equilibrium of an elastic body in frictional or frictionless contact with a foundation. Here we restrict ourselves to present only one example and, to keep this paper
in a reasonable length, we restrict to the homogeneous case, skip the description of the model, and refer the reader to the books \cite{C, DL, KO, SofMat, SMBOOK} for details on the physical setting of contact problems, statement of the models and various mechanical interpretation.

To introduce the  problem  we need the following notations. First,  $\mathbb{S}^d$ will represent the space of second order symmetric tensors on $\mathbb{R}^d$ and $``\cdot"$,  $\|\cdot\|$, $\bzero$ will denote the canonical inner product, the Euclidian norm and the zero element of  $\mathbb{R}^d$  and $\mathbb{S}^d$, respectively. Let
$\Omega\subset\mathbb{R}^d$ ($d=2,3$) be a smooth domain with outward normal $\bnu$ and boundary $\Gamma$ and let
$\Gamma_1$, $\Gamma_2$, $\Gamma_3$ be a partition of $\Gamma$ such that ${ meas}\,(\Gamma_1)>0$. Denote
\begin{eqnarray*}
	&&X=\{\,\bv\in H^1(\Omega)^d:\  \bv =\bzero\ \ {\rm on}\ \Gamma_1\,\}.
\end{eqnarray*}
It is well known that  $X$ is a Hilbert space 
endowed with the inner product
\begin{equation}
(\bu,\bv)_X= \int_{\Omega}
\bvarepsilon(\bu)\cdot\bvarepsilon(\bv)\,dx,
\end{equation}
where $\bvarepsilon(\bv)$ denotes the linearized strain of  $\bv$, i.e. $\bvarepsilon(\bv)=\frac{1}{2}(\nabla\bv+\nabla^T\bv)$.
We use $\bzero_X$ for the zero element of $X$ and, for any element $\bv\in X$,  we still write $\bv$ for the trace of $\bv$ to $\Gamma$. Moreover, we
denote by $v_\nu$ and $\bv_\tau$ its normal and
tangential components given by
$v_\nu=\bv\cdot\bnu$ and $\bv_\tau=\bv-v_\nu\bnu$, respectively, 
Finally, we recall that  the Sobolev
trace theorem yields
\begin{equation}\label{trace}
\|\bv\|_{L^2(\Gamma_3)^d}\le \| \gamma \|\,\|\bv\|_{X}\quad
\forall\,\bv \in X.
\end{equation}
where, here and below, $\| \gamma \|$ denotes the norm of the trace
operator $\gamma \colon X \to L^2(\Gamma_3)^d$.

\medskip
Consider in what follows the data ${\cal F}$, $B$, $F$, 
$p$,  $k$, $\fb_0$, $\fb_2$, $\omega$ and $g$, assumed to satisfy the following conditions.

\begin{equation}
\left\{\begin{array}{ll} {\cal F}\colon 
\mathbb{S}^d\to \mathbb{S}^d\ \mbox{is such that} \\ [1mm]
{\rm (a)\  there\ exists}\ L_{\cal F}>0\ {\rm such\ that}\\
{}\qquad \|{\cal F}\bvarepsilon_1-{\cal F}\bvarepsilon_2\|
\le L_{\cal F} \|\bvarepsilon_1-\bvarepsilon_2\|\\
{}\qquad \quad\mbox{for all} \ \ \bvarepsilon_1,\bvarepsilon_2
\in \mathbb{S}^d; 
\\ [1mm]
{\rm (b)\  there\ exists}\ m_{\cal F}>0\ {\rm such\ that}\\
{}\qquad ({\cal F}\bvarepsilon_1-{\cal F}\bvarepsilon_2)
\cdot(\bvarepsilon_1-\bvarepsilon_2)\ge m_{\cal F}\,
\|\bvarepsilon_1-\bvarepsilon_2\|^2\quad\\
{}\qquad\quad \mbox{for all} \ \ \bvarepsilon_1,
\bvarepsilon_2 \in \mathbb{S}^d.\\ [1mm]
\end{array}\right.
\label{m1}
\end{equation}

\begin{eqnarray}
&&\label{m2} B \ \ \mbox{is a closed  convex subset of} \ \mathbb{S}^d\ \mbox{such that}\ \bzero\in B.
\end{eqnarray}

\begin{equation}
\left\{\begin{array}{ll} { F}\colon 
\mathbb{R}\to \mathbb{R}\ \mbox{is such that} \\ [1mm]
{\rm (a)\  there\ exists}\ L_{F}>0\ {\rm such\ that}\\
{}\qquad |F(r_1)-F(r_2)|
\le L_{F} |r_1-r_2|\quad\mbox{for all} \ r_1,r_2
\in \mathbb{R}; 
\\ [1mm]
{\rm (b)\  }\ F(r)=0\quad\mbox{for all}\ r\le 0.
\end{array}\right.
\label{m3}
\end{equation}

\medskip
\begin{equation}
\left\{\begin{array}{ll} {p}\colon 
\mathbb{R}\to \mathbb{R}\ \mbox{is such that} \\ [1mm]
{\rm (a)\  there\ exists}\ L_{p}>0\ {\rm such\ that}\\
{}\qquad |p(r_1)-p(r_2)|
\le L_{p} |r_1-r_2|\quad\mbox{for all} \ r_1,r_2
\in \mathbb{R}; \ \
\\ [1mm]
{\rm (b)\  }\ p(r)=0\quad\mbox{for all}\ r\le 0.
\end{array}\right.
\label{m4}
\end{equation}

\begin{eqnarray}
&&\label{m4m}(L_F+L_p)\| \gamma \|^2<m_{\cal F}.\\ [2mm]
&&\label{m5}k\in L^2(\Gamma_3),\qquad k(\bx)\ge 0\quad{\rm a.e.} \ \bx\in\Gamma_3.\\[2mm]
&&\label{m6} \fb_0 \in L^2(\Omega)^d, \quad\  \fb_2\in L^2(\Gamma_2)^d. \\[2mm]
&&\label{m7}\omega\in L^\infty(\Omega),\qquad \omega(\bx)\ge 0\quad{\rm a.e.} \ \bx\in\Omega.\\[2mm]
&&\label{m8}g\in L^2(\Gamma_3),\qquad 0\le g(\bx)\le k(\bx)\quad{\rm a.e.} \ \bx\in\Gamma_3.
\end{eqnarray}

\medskip
We denote by $P_B:\mathbb{S}^d\to B$ the projection operator on $B$ and we consider the function $q:\real\to\real$ defined by
\begin{equation}\label{jr}
q(r)=\int_0^r p(s)\,ds\qquad\forall\,r\in\real.
\end{equation}
Note that (\ref{m4}) shows that the function $p$ is Lipschitz continuous. Nevertheless, it could be nonmonotone and, as a result the function $q$ could be nonconvex. 
\medskip

Let $K$, $A$, $\varphi$, $j$, $\fb$ be defined as follows:
\begin{eqnarray}
\label{8b0}&&K=\{\,\bv\in X\ :\ v_\nu \le k\ \  \hbox{a.e. on}\
\Gamma_3\,\},\\ [3mm]
&&A\colon X \to X^*,\quad
\label{8b1}\langle A\bu,\bv\rangle =\int_{\Omega}\cF\bvarepsilon(\bu)\cdot\bvarepsilon(\bv)\,dx\\ [2mm]
&&\qquad\qquad\qquad+\int_{\Omega}\omega\big(\bvarepsilon(\bu)-P_B\bvarepsilon(\bu)\big)\cdot\bvarepsilon(\bv)\,dx,\nonumber\\[2mm]
&&\varphi\colon X\times X \to \real, 
\label{8b3}\quad
\varphi(\bu,\bv)=\int_{\Gamma_3} F(u_\nu-g)\,\|\bv_\tau\|\,da,\\[2mm]
&&j\colon X\to\real, \quad
\label{8b5}j(\bv)=\int_{\Gamma_3}q(v_\nu-g)\, da,\\[2mm]
&&\fb\in X^*,\quad\label{8ef}\langle\fb,\bv\rangle
=\int_{\Omega}\fb_0\cdot\bv\,dx +
\int_{\Gamma_2}\fb_2\cdot\bv\,da,
\end{eqnarray}
for all $\bu,\bv\in X$.
With these notations we consider the following
problem.

\medskip\noindent{\bf Problem} ${\cal P}$. {\it	Find a function $\bu\in K$ such that}
\begin{eqnarray}
&&\label{8hv}
\langle A\bu,\bv - \bu \rangle + \varphi (\bu,\bv) - \varphi (\bu,\bu) + j^0(\bu;\bv -\bu)
\\[2mm]
&&\qquad\qquad\ge \langle \fb, \bv - \bu \rangle 
\ \ \mbox{for all} \ \ \bv \in K.\nonumber
\end{eqnarray}

Problem ${\cal P}$  represents the variational formulation of a mathematical model which describes the equilibrium of an elastic body  in frictional contact with a foundation made of a rigid body covered by a layer of deformable material, say asperities. The body is fixed on $\Gamma_1$, is acted by body forces and surface tractions on $\Gamma_2$, and is in potential contact on $\Gamma_3$ with a foundation.
The functions  ${\cal F}$, $\omega$ and the set $B$ are related to the constitutive law of the material, $\fb_0$ and $\fb_2$ denote the density of body forces and  surface tractions, respectively, $g$ is the initial gap, and $k-g$ represents the thickness of the deformable material. The function $p$ is the so-called normal compliance function which  describes the behaviour of the deformable layer of the foundation and $F$ represents the friction bound. Note that part of the assumptions on these data presented above are not necessary from mathematical point of view. However, we adopt them since they are imposed for mechanical reasons.

A problem similar to Problem ${\cal P}$ was considered in \cite[Chapter 7]{SMBOOK}, with $g\equiv\omega\equiv0$. Nevertheless, there, the nonhomogeneous case was considered, i.e., the  functions ${\cal F}$, $F$ and $p$ were supposed to depend on the spatial variable $\bx$. Since the case when $\omega$ and $g$ are  positive functions does not introduce important modification, we skip the proof of  the following result and send the reader to \cite{SMBOOK} for more details.

\begin{Lemma}\label{l30}
Assume $(\ref{m1})$--$(\ref{m8})$. Then the set $K$, operator $A$, and functions $\varphi$ and $j$ satisfy assumptions  $(\ref{K})$--$(\ref{smal})$.
\end{Lemma}

We now illustrate the use of the abstract results in Theorems~\ref{t1} and \ref{t3}
in the study of Problem~ ${\cal P}$.

\begin{Theorem}\label{t6}
	Assume  $(\ref{m1})$--$(\ref{m8})$.
	Then the following statements hold.
	
\medskip	
{\rm  a)} 	Problem~${\cal P}$ has a unique solution $\bu \in  K$.

\medskip
{\rm  b)}  Problem~${\cal P}$  is well-posed.

\medskip
{\rm  c)} The solution of Problem  ${\cal P}$  depends continuously on $\fb_0$, $\fb_2$ and $g$, i.e., if $\bu_n$ represents the solution of Problem ${\cal P}$ with the data  $\fb_{0n}$, $\fb_{2n}$ and $g_n$ which have the regularity prescribed in $(\ref{m6})$, $(\ref{m8})$ and
\begin{eqnarray}
&&\label{cw1}\fb_{0n}\to\fb_0\quad{\rm in}\quad L^2(\Omega)^d,\quad \fb_{2n}\to\fb_2\quad{\rm in}\quad L^2(\Gamma_2)^d,\\[2mm]
&&\label{cw2}g_n\to g\quad{\rm in}\quad L^2(\Gamma_3),
\end{eqnarray}
then
\begin{equation}
\label{cw} \bu_n\to \bu\quad\mbox{\rm in}\quad  X. 
\end{equation}
\end{Theorem}

\medskip

\medskip\noindent{\it Proof.}  Part ${\rm  
	a)}$  is a direct consequence Theorem \ref{t5} and Lemma \ref{l30}. Moreover, part  ${\rm  b)}$ follows from Corollary \ref{c2}.

For  part ${\rm  c)}$ we use Theorem \ref{t3}.
To this end we assume in what follows that
$\fb_{0n}$, $\fb_{2n}$, and $g_n$  satisfy $(\ref{m6})$ and $(\ref{m8})$ and, for each $n\in\mathbb{N}$, we consider  the functions $\varphi_n$, $j_n$, and the element $\fb_n\in X^*$ defined by 
\begin{eqnarray}
&&\varphi_n\colon X\times X \to \real, 
\label{8b3n}\quad
\varphi_n(\bu,\bv)=\int_{\Gamma_3} F(u_\nu-g_n)\,\|\bv_\tau\|\,da,\\[2mm]
&&j_n\colon X\to\real, \quad
\label{8b5n}j_n(\bv)=\int_{\Gamma_3}q(v_\nu-g_n)\, da,\\[2mm]
&&\fb_n\in X^*,\quad\label{8efn}\langle\fb_n,\bv\rangle
=\int_{\Omega}\fb_{0n}\cdot\bv\,dx +
\int_{\Gamma_2}\fb_{2n}\cdot\bv\,da,
\end{eqnarray}
for all $\bu,\bv\in X$. 

Let $\bu,\, \bv\in X$. Then, using assumption (\ref{m3}) on the function $F$  it follows that
\begin{eqnarray*}
&&\varphi_n(\bu,\bv)-\varphi_n(\bu,\bu)-
\varphi(\bu,\bv)+\varphi(\bu,\bu)\\ [3mm]
&&\quad=\int_{\Gamma_3}\big( F(u_\nu-g_n)-F(u_\nu-g)\big)\big(\|\bv_\tau\|-\|\bu_\tau\|\big)\,da
\\ [2mm]
&&\qquad\le L_F\int_{\Gamma_3}|g_n-g|\,\|\bu_\tau-\bv_\tau\|\,da
\end{eqnarray*}
and, therefore, the trace inequality (\ref{trace}) yields
\begin{eqnarray}
	&&\label{n2}\varphi_n(\bu,\bv)-\varphi_n(\bu,\bu)-
	\varphi(\bu,\bv)+\varphi(\bu,\bu)\\ [3mm]
	&&\quad\le L_F\|\gamma\|\,\|g_n-g\|_{L^2(\Gamma_3)}\|\bu-\bv\|_X. \nonumber
\end{eqnarray}
Next, using (\ref{m4}) and (\ref{jr}) it follows that the functions $j_n$ and $j$ are regular, and, moreover,
\[j_n^0(\bu ;\bv-\bu)=\int_{\Gamma_3}p(u_\nu-g_n)(v_\nu-u_\nu)\,da,\quad j^0(\bu;\bv-\bu)=\int_{\Gamma_3}p(u_\nu-g)(v_\nu-u_\nu)\,da.\]
Therefore, using arguments similar to those used in the proof for (\ref{n2}) we deduce that
\begin{equation}\label{n3}
	j_n^0(\bu;\bv-\bu)-j^0(\bu;\bv-\bu)\le L_p\|\gamma\|\,\|g_n-g\|_{L^2(\Gamma_3)} \|\bu-\bv\|_X.
\end{equation}
It follows from (\ref{n2}) and (\ref{n3}) that conditions
(\ref{fin}), (\ref{jn}) hold with
\[
b_n=L_F\|\gamma\|\,\|g_n-g\|_{L^2(\Gamma_3)},\quad
c_n=L_p\|\gamma\|\,\|g_n-g\|_{L^2(\Gamma_3)},
\]
and, using assumptions (\ref{cw2}), we deduce that conditions (\ref{cn}) and (\ref{ccn}) are satisfied.
On the other hand, it is easy to see that the convergences (\ref{cw1}) imply (\ref{cfn}) for $\fb_n$ and $\fb$ given by
(\ref{8efn}) and (\ref{8ef}), respectively.
Finally, recall that part {\rm (i)} of the theorem gurantees that the variational-hemivariational inequality (\ref{8hv}) is 
well-posed. We are now in a position to apply Theorem \ref{t3}  in order to deduce the convergence (\ref{cw}) which concludes the proof. \hfill$\Box$

\medskip

In addition to the mathematical interest in the convergence result
in Theorem~\ref{t3} c), it is important from the mechanical point of view, since it provides the  continuous dependence of the solution with respect to the density of the body forces and tractions and the gap function.

\medskip
We end this section with the remark that the strongly monotonicity of the operator $A$, guaranteed by condition (\ref{m1}), plays a crucial role in the well-posedness of Problem ${\cal P}$. If this condition does not hold, in general, Problem ${\cal P}$ is not well-posed. To provide an example, assume in what follows that ${\cal F}$ vanishes, and consider the set $\widetilde{K}\subset K$ defined by
\[\widetilde{K}=\{\,\bv\in K\ : \bvarepsilon(\bu)\in B \  \ \hbox{a.e. in}\ \Omega,\  \ v_\nu \le g\ \  \hbox{a.e. on}\ \Gamma_3\,\}\]
and assume that the body forces and tractions vanish, which implies that $\fb=\bzero_{X^*}$.
It is easy to see that in this case $A\bu=\bzero_{X^*}$, $\varphi(\bu,\bv)=0$, $j^0(\bu,\bv)=0$ for all $\bu\in \widetilde{K}$, $\bv\in X$ and, therefore, any element $\bu\in \widetilde{K}$ is a solution ot Problem ${\cal P}$. Moreover, it follows from assumption (\ref{m2}) that $\bzero_X\in\widetilde{K}$. On the other hand, concrete examples of convex sets $B$ and reference configurations $\Omega$ for which $\widetilde{K}$ contains at least one element $\bu\ne\bzero_X$ can be easily provided. We deduce from here that
in this case Problem ${\cal P}$
has more than one solution  and, therefore, is not well-posed.

\medskip

\section*{Acknowledgments}

\indent This project has received funding from the European Union's Horizon 2020
Research and Innovation Programme under the Marie Sklodowska-Curie
Grant Agreement No 823731 CONMECH.
This research was also supported by the National Natural Science Foundation of China (11771067) and the Applied Basic Project of Sichuan Province (2019YJ0204).

\end{document}